\documentclass[reqno,11pt]{amsart}

\usepackage{epsf}
\usepackage{graphics}
\usepackage{graphicx}
\usepackage{amssymb}
\usepackage{amsmath}
\usepackage{mathcomp}
\usepackage{wasysym}
\date{}

\usepackage{microtype,hyperref}
\hypersetup{colorlinks={true},linkcolor={black},citecolor={black},filecolor={black},urlcolor={black},pdfauthor={\authors},pdftitle={\title}}

\theoremstyle{plain}
\newtheorem{theorem}{Theorem}

\newtheorem{lemma}{Lemma}

\theoremstyle{definition}

\theoremstyle{remark}

\def\N{{\mathbb N}}
\def\Z{{\mathbb Z}}
\def\Q{{\mathbb Q}}

\title{Average four-genus of two-bridge knots}

\author{S.~Baader, A.~Kjuchukova, L.~Lewark, F.~Misev, A.~Ray}

\begin{document}

\begin{abstract} We prove that the expected value of the ratio between the smooth four-genus and the Seifert genus of two-bridge knots tends to zero as the crossing number tends to infinity. 
\end{abstract}

\maketitle

\section{Introduction}

The 4-genus $g_4(K)$ of a knot~$K \subset S^3$ is defined to be the minimal genus among all smoothly embedded surfaces $\Sigma \subset D^4$ with boundary $\partial \Sigma=K$. Unlike the classical Seifert genus $g(K)$, the 4-genus can be zero for non-trivial knots. In particular, the 4-genus can be strictly smaller than the Seifert genus. Attempts to  determine the 4-genus for special classes of knots tend to fail. A famous counterexample to this rule is the class of torus knots, for which the 4-genus was determined by elaborate methods (see~\cite{KM,Ra}, and~\cite{Ru} for the more general case of quasipositive knots). Even the classification of slice knots, i.e.\ knots with vanishing 4-genus, turns out to be a very difficult task and has only been achieved for special classes of knots, such as 2-bridge knots and certain types of Montesinos knots (see~\cite{Li} and~\cite{Le}, \cite{GJ}).

In this paper, we determine the expected value of the ratio of the 4-genus and the Seifert genus for large 2-bridge knots. More precisely, for all natural numbers~$n \geq 2$, we define $\left \langle \frac{g_4}{g} \right \rangle_n \in \Q$ to be the average value of the ratio $g_4(K)/g(K)$, taken over all 2-bridge knots~$K$ whose 4-plat presentation, in the sense of Conway~\cite{Co}, has~$2n$ crossings.

\begin{theorem} 
$$\lim_{n \to \infty} \left \langle \frac{g_4}{g} \right \rangle_n=0.$$
\end{theorem}

It would be interesting to consider other natural classes of knots, such as alternating knots or all knots, for which we expect similar results to hold.
However, one needs to restrict to classes of links that admit arbitrarily small ratios $g_4/g$. In the extreme case of strongly quasipositive knots, this ratio is always one~\cite{Ru}. Throughout this paper, we use the following notation for a 2-bridge knot in Conway's 4-plat presentation~\cite{Co}. Let $m \in \N$ be a natural number and let $a_1,a_2,\ldots, a_{2m} \in\Z$ be non-zero integers. We define the oriented 2-bridge knot $K(2a_1,\ldots,2a_{2m})$ via the standard diagram $C(2a_1,\ldots,2a_{2m})$ with~$2m$ twist regions, as shown in Figure~1.
\begin{figure}[ht]
\scalebox{1.0}{\raisebox{-0pt}{$\vcenter{\hbox{\includegraphics{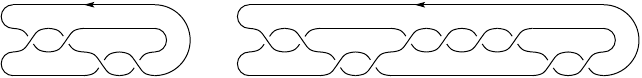}}}$}}
\caption{2-bridge knot diagrams $C(2,2)$ and $C(-2,-2,4,-2)$}
\end{figure}

The sign convention for crossings is such that positive parameters $a_1,a_3,\dots$ and negative parameters $a_2,a_4,\dots$ stand for positive crossings. As a result, the diagram $C(2a_1,\ldots,2a_{2m})$ is alternating if and only if all the parameters have the same sign. The simplest non-trivial 2-bridge knots are the figure-eight knot $K(2,2)$ and the trefoil knot $K(2,-2)$. All 2-bridge knots admit a 4-plat presentation with an even number of crossings in each twist region, a fact that is easy to prove in the framework of finite continued fractions (see Chapter~2 in~\cite{Ka}). Restricting to even parameters has two major advantages:
\begin{enumerate}
\item the Seifert genus satisfies $g(K(2a_1,\ldots,2a_{2m}))=m$,
\item the representation $K(2a_1,\ldots,2a_{2m})$ is unique, up to the symmetry
$K(2a_1,\ldots,2a_{2m})=K(-2a_{2m},\ldots,-2a_1)$. 
\end{enumerate}
The first fact follows from the Alexander polynomial degree bound, or more geometrically from Gabai's result on the plumbing structure of the canonical Seifert surface associated with the diagram $C(2a_1,\ldots,2a_{2m})$~\cite{Ga}. The second fact follows again from properties of continued fractions, in combination with the original classification of 2-bridge knots via branched coverings (see~\cite{Re,Sch}, as well as~\cite{Ka} and the much more general classification of arborescent knots~\cite{BS}).

Recall that any $2$-bridge knot is alternating. The Conway 4-plat presentation does not in general minimise the crossing number, unless it is itself alternating. Nevertheless, the number $2(|a_1|+|a_2|+\ldots+|a_{2m}|)$ is a fairly accurate approximation of the crossing number of the knots $K(2a_1,\ldots,2a_{2m})$. In the worst case, there are non-alternating transitions between all the $2m$ twist regions. Each of these transitions can be turned into an alternating one by a move that removes one crossing, as sketched in Figure~2. Therefore, every 4-plat diagram can be transformed into a reduced alternating diagram by removing at most $2m-1$ crossings. In particular, the crossing number $c(K(2a_1,\ldots,2a_{2m}))$ is larger than $|a_1|+|a_2|+\ldots+|a_{2m}|$. With this in mind, we assert without proof that Theorem~1 holds with $n$ replaced by the crossing number, i.e.\ by taking the average value of the ratio $g_4(K)/g(K)$ over all 2-bridge knots~$K$ with crossing number~$n$.

\begin{figure}[ht]
\scalebox{1.0}{\raisebox{-0pt}{$\vcenter{\hbox{\includegraphics{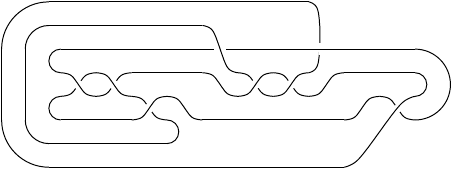}}}$}}
\caption{Reducing the diagram $C(-2,-2,4,-2)$}
\end{figure}

In order to estimate the average ratio $\left \langle \frac{g_4}{g} \right \rangle_n$, we will consider $\langle g \rangle_n \in \Q$ and $\langle g_4 \rangle_n \in \Q$, the average values of the Seifert genus and the 4-genus respectively, taken over all 2-bridge knots $K(2a_1,\ldots,2a_{2m})$ with $n=|a_1|+|a_2|+\ldots+|a_{2m}| \geq 2$. The following two statements almost immediately imply Theorem~1, as we will see in the next section.

\begin{lemma}\label{lem-1}
$$\liminf_{n \to \infty} \frac{\langle g \rangle_n}{n} \geq \frac{1}{4}.$$
\end{lemma}

\begin{lemma}\label{lem-2}
$$\lim_{n \to \infty} \frac{\langle g_4 \rangle_n}{n}=0.$$
\end{lemma}

\section{Proof of Lemma~\ref{lem-1}: Average Seifert genus}

Facts~(1) and~(2) in the previous section allow for an easy estimate of the average genus of 2-bridge knots. Indeed, up to signs of all the parameters $a_k$ and up to the symmetry $K(2a_1,\ldots,2a_{2m})=K(-2a_{2m},\ldots,-2a_1)$, the number of 2-bridge knots with fixed $n=|a_1|+|a_2|+\ldots+|a_{2m}| \geq 2$ is equal to the number of ordered partitions of the natural number $n$ into an even number of strictly positive natural numbers. For fixed $m \in \N$, there are ${n-1} \choose {2m-1}$ ordered partitions of $n$ into precisely~$2m$ strictly positive natural numbers. Since we are interested in the average value of the classical genus $g(K(2a_1,\ldots,2a_{2m}))=m$, we need to determine the average value of $m$ over all ordered even partitions of $n$, as $n$ tends to infinity. The result is seen to be $(n+1)/4$, using the symmetry ${{n-1} \choose {2m-1}}={{n-1} \choose {n-2m}}$ and solving for the middle coefficient: $2m-1=\frac{n-1}{2}$. Here we are a bit imprecise since we disregarded the signs of the parameters $a_k$ (as well as the symmetry~(2); however, that does not affect the average value of~$m$). Taking these signs into account causes a shift of the average value of~$m$ in the positive direction, since the number of signed ordered partitions of~$n$ into $2m$ numbers is now multiplied by $2^{2m}$. In particular, there are more signed partitions with $m=\frac{n+1}{4}+k$ than signed partitions with $m=\frac{n+1}{4}-k$, so the average value is at least~$\frac{n+1}{4}$. This settles Lemma~\ref{lem-1}, and even more. The peak of the distribution of the binomial coefficients gets narrower as $n$ increases. More precisely, the proportion of ordered partitions of $n$ with $m \leq n/8$ tends to zero as $n$ tends to infinity. This is a weak consequence of the local limit theorem of de Moivre-Laplace (see Chapter~3.1 in~\cite{Si}). As a consequence, we can neglect all 2-bridge knots with Seifert genus smaller than $n/8$ in all asymptotic considerations. Here we are using Fact (1). This together with Lemma~\ref{lem-2} implies Theorem~1:
$$\lim_{n \to \infty} \left\langle \frac{g_4}{g} \right\rangle_n \leq \lim_{n \to \infty}\frac{8\langle g_4 \rangle_n}{n}=0.$$

\section{Proof of Lemma~\ref{lem-2}: Average 4-genus}

In this section, we prove that the average value of the 4-genus of 2-bridge knots has sublinear growth in $n=|a_1|+\ldots+|a_{2m}|$. For this purpose we need the following two facts:
\begin{enumerate}
\item[(i)] for all knots~$K$, the connected sum of~$K$ with its mirror image $\bar{K}$ is a ribbon knot; in particular $g_4(K \# \bar{K})=0$,
\item[(ii)] for all $i \leq m$, the knot $K(2a_1,\ldots,2a_{2m})$ is related to the connected sum $K(2a_1,\ldots,2a_{2i-2}) \# K(2a_{2i+1},\ldots,2a_{2m})$ by an oriented cobordism of genus~$1$, as shown in Figure~3.
\end{enumerate}

\begin{figure}[ht]
\scalebox{1.0}{\raisebox{-0pt}{$\vcenter{\hbox{\includegraphics{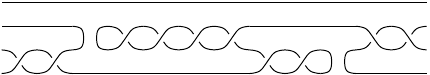}}}$}}
\caption{Cobordism of genus one (two saddle moves)}
\end{figure}

In the preceding section, we saw that the average value of~$m$, taken over all ordered partitions of~$n$ with~$2m$ strictly positive natural numbers, is~$(n+1)/4$. Moreover, the proportion of partitions with $m<n/8$ tends to zero as~$n$ tends to infinity. Therefore, in the process of determining the average 4-genus of 2-bridge knots with large~$n$, we may again ignore the contribution of all knots with $m<n/8$, i.e.\ of all knots with fewer than $n/8$ twist regions. From now on, we consider 2-bridge knots $K=K(2a_1,\ldots,2a_{2m})$ with $m \geq n/8$ only.

Let $k \geq 32$ be a constant smaller than~$n$. We remove all pairs of consecutive twist regions where at least one of $|a_{2i-1}|$ or $|a_{2i}|$ is larger than~$k$, by an oriented cobordism of genus~$2$. For this, we first apply the cobordism shown in Figure~3, untwist the two ribbons in question, then apply the inverse of that cobordism. There are at most $n/k$ pairs of twist regions of that type, since $n=|a_1|+|a_2|+\ldots+|a_{2m}|$. An oriented cobordism of genus at most $2n/k$ relates our initial knot $K$ with a new 2-bridge knot $K'$, with at least $n/8-2n/k \geq n/16$ twist regions and all parameters $a_j \leq k$. 

Now we fix another constant $s \geq 2$ smaller than~$n$. Repeated use of (ii) at all multiples of~$s+1$, i.e.\ for $i=s+1,2(s+1),3(s+1),\ldots$ transforms $K'$ into a connected sum of 2-bridge knots $K_1 \# K_2 \#\ldots \# K_t$. Each summand has precisely $2s$ twist regions, except possibly $K_t$, which might be shorter, but that does not matter for our consideration. The number $t$ is the ceiling of $m/(s+1)$, which implies, up to an error of 1,
$$\frac{n}{16(s+1)} \leq t \leq \frac{n}{2(s+1)},$$
since the number $m$ of twist regions satisfies $n/16 \leq m \leq n/2$. A cobordism of genus~$t \leq \frac{n}{2(s+1)}$ relates the knot $K'$ and the connected sum of knots $K''=K_1 \# K_2 \#\ldots \# K_t$. There are $(2k)^{2s}$ types of 2-bridge knots with $2s$ twist regions and all parameters bounded in absolute value by~$k$, since these can be seen as words of length~$2s$ in the alphabet $\{-k,\ldots,-1,1,\ldots,k\}$. Fix a word $w=b_1 b_2 \cdots b_{2s}$ in that alphabet, and let $-w=(-b_1) (-b_2) \cdots (-b_{2s})$ denote its mirror image. The knot $K''$ contains a certain number $a(w)$ of summands of type $K(w)$, and a number $a(-w)$ of summands of type $K(-w)$. The expected value of $a(w)$, as well as that of $a(-w)$, is $\frac{t}{(2k)^{2s}}$. Now comes the core part of the argument:

The expected value of the difference $|a(w)-a(-w)|$ is $\sqrt{\frac{t}{(2k)^{2s}}}$, which is at most $\sqrt{\frac{n}{2(s+1)(2k)^{2s}}}$. This fact is a well-known feature of random walks on~$\Z$: the expected distance between the initial and final position of a random walk on $\Z$ after $\beth$ steps of size $1$ is $\sqrt{\beth}$. This is another consequence of the local limit theorem: the binomial distribution ${{m} \choose {k}}2^{-m}$ is well-approximated by 
$$\sqrt{\frac{2}{\pi m}} e^{-\frac{2}{m}(k-\frac{m}{2})^2},$$
in an interval of size proportional to $\sqrt{m}$ around $\frac{m}{2}$, and negligible in the complement (see Chapter~3.1 in~\cite{Si} for a precise statement). In our situation, this means that the knot $K''$ contains an average number of $\frac{t}{(2k)^{2s}}$ copies of the ribbon knot $K(w)\#K(-w)$, plus a few extra summands of the form $K(w)$ or $K(-w)$, but only at most $\sqrt{\frac{n}{2(s+1)(2k)^{2s}}}$ many. Formally, we could interpret the sequence of knots $K_1,K_2,\ldots,K_t$ as a random walk of length~$t$ on~$\Z^{(2k)^{2s}/2}$, where each coordinate direction stands for a pair of words $\{w,-w\}$. The expected value of the $L^1$-norm of the endpoint is bounded above by $\sqrt{t}$. This means that we have at most $\sqrt{t} \leq \sqrt{\frac{n}{2(s+1)}}$ extra summands to remove. As a consequence, we obtain an oriented cobordism of (expected) genus at most $s \sqrt{\frac{n}{2(s+1)}}$ between the knot $K''$ and a ribbon knot. Altogether, we obtain an oriented cobordism of genus at most
$$\frac{2n}{k}+\frac{n}{2(s+1)}+s \sqrt{\frac{n}{2(s+1)}}$$
between the initial knot $K$ and a ribbon knot. The expected value of the 4-genus of $K$ is therefore bounded by that same quantity. Dividing by $n$ and letting $n$ tend to infinity yields $\frac{2}{k}+\frac{1}{2(s+1)}$.  The constants~$k$ and~$s$ can be chosen arbitrarily large, thus we get the desired limit:
$$\lim_{n \to \infty} \frac{\langle g_4 \rangle_n}{n}=0.$$

\subsection*{Acknowledgements}
Part of this work was completed when SB visited AK, FM, and AR at the Max Planck Institute for Mathematics. We are grateful to the MPIM for its support and hospitality.

\bigskip
\noindent
Mathematisches Institut, Sidlerstr.~5, 3012 Bern, Switzerland

\smallskip
\noindent
Max-Planck-Institut f\"ur Mathematik, Vivatsgasse 7, 53111 Bonn, Germany

\smallskip
\noindent
Mathematisches Institut, Sidlerstr.~5, 3012 Bern, Switzerland

\smallskip
\noindent
Max-Planck-Institut f\"ur Mathematik, Vivatsgasse 7, 53111 Bonn, Germany

\smallskip
\noindent
Max-Planck-Institut f\"ur Mathematik, Vivatsgasse 7, 53111 Bonn, Germany

\newcommand{\myemail}[1]{\texttt{\href{mailto:#1}{#1}}}

\bigskip

\smallskip
\noindent
\myemail{sebastian.baader@math.unibe.ch}

\smallskip
\noindent
\myemail{sashka@mpim-bonn.mpg.de}

\smallskip
\noindent
\myemail{lukas@lewark.de}

\smallskip
\noindent
\myemail{fmisev@mpim-bonn.mpg.de}

\smallskip
\noindent
\myemail{aruray@mpim-bonn.mpg.de}

\end{document}